\input amstex.tex 
\magnification=\magstep{1.5}
\baselineskip 22pt
\documentstyle{amsppt}
\topmatter
\title 
Automorphism groups in a family of K3 surfaces 
\endtitle
\author 
Keiji Oguiso 
\endauthor 
\address 
Department of Mathematical Sciences University of Tokyo
153-8914 Komaba Meguro Tokyo Japan
\endaddress
\email
oguiso\@ms.u-tokyo.ac.jp
\endemail
\subjclass 
14J28 
\endsubjclass
\abstract 
A few facts concerning the phrase "the automorphism groups become 
larger at special points of the moduli of K3 surfaces" are presented. 
It is also shown that the automorphism groups are of infinite order over a dense subset in any one-dimensional non-trivial family of projective K3 surfaces.  
\endabstract
\leftheadtext{Keiji Oguiso}
\rightheadtext{Automorphism groups of K3 surfaces}
\endtopmatter
\document
\head
{Introduction}
\endhead 
According to the global Torelli Theorem for K3 surfaces due to 
Piatetski-Shapiro and Shafarevich [PSS] (See also [BPV]), one can describe 
the automorphism group of a K3 surface as the quotient of the orthogonal 
group of the Picard 
lattice by the $2$-refection group, up to finite groups, i.e. up to finite 
kernel and cokernel, which we denote by $\equiv$ later. Although an explicit 
calculation of such a quotient group is 
rather hard even in a concrete case, beautiful calculations of the 
automorphism groups of some special but quite remarkable K3 surfaces have been 
done by Vinberg, Keum, Kondo and others: for instance those of K3 surfaces of $\rho = 20$ and of smaller discriminant by [Vi] and [KK], those of generic Kummer surfaces by [Ke] and 
[Kn3]. In these examples, the automorphism groups are of infinite order. 
As a counter part, Kondo [Kn1] has also classified the 
automorphism groups of K3 surfaces $X$ with $\vert \text{Aut}(X) \vert 
< \infty$ based on Nukulin's classification of both the Picard lattices and the dual graphs of smooth rational curves, i.e. generators and relations of the 
$2$-reflection groups of such K3 surfaces ([Ni5] and references therein). 
\par 
\vskip 4pt 
On the other hand, since no polarization on a K3 surface is invariant 
under $\text{Aut}(X)$ unless $\vert \text{Aut}(X) \vert < \infty$ and since 
the behaviour of the Picard lattices is rather unstable under deformation (see 
for instance [Og]), even when one considers a projective family, one can not 
relate the full automorphism groups $\text{Aut}(X)$ well with the 
automorphism group $\text{PGL}(n)$ of the ambient space $\bold P^{n}$ nor 
simultaneously with 
the orthogonal groups of their Picard lattices, and it does 
not seem so clear 
how $\text{Aut}(X)$ behaves when a K3 surface $X$ varies in a family. 
\par 
\vskip 4pt 
Aim of this short note is to clarify certain behaviour of the automorphism 
groups of K3 surfaces under deformation. Saying more explicitly, our interest is to 
discuss about the phrase, "{\it the 
automorphism groups of K3 surfaces become larger at special points in their 
moduli.}" (See the main Theorem and Corollary below.)
\par
\vskip 4pt
Here and hereafter, we shall work over the complex number field $\bold C$. 
By a K3 surface we mean a simply-connected smooth {\it projective} surface $X$ 
which admits an everywhere non-vanishing holomorphic $2$-from $\omega_{X}$. 
\par 
\vskip 4pt 
In order to state the results in the strongest form, throughout this note, we 
shall consider a non-trivial smooth family of K3 surfaces 
over the unit disk: 
$$\varphi : \Cal X \rightarrow \Delta := \{ t \in \bold C \vert\,\, 
\vert t \vert 
< 1 \}.$$ 
If prefer, one may regard this family either as a 
$1$-dimensional small part of the Kuranishi family of a given K3 surface or 
as a $1$-dimensional small non-trivial part of a global family of K3 surfaces, 
e.g. a universal family of polarized K3 surfaces of given degree. We regard $\Delta$ as a topological space by the 
Euclidean topology unless stated otherwise. We denote the fiber 
$\varphi^{-1}(t)$ by $\Cal X_{t}$. 
\par 
\vskip 4pt 
For the statement, we choose a marking
$$\tau : R^{2}\varphi_{*}\bold Z_{\Cal X} \rightarrow \Lambda \times 
\Delta.$$ 
Here $\Lambda$ is the K3 lattice, i.e. the lattice isomorphic to 
$U^{\oplus 3} \oplus E_{8}(-1)^{\oplus 2}$. 
Set $\Lambda_{t} := \tau_{t}(\text{NS}(\Cal X_{t}))$. By [Og], 
one can find a primitive sublattice $\Lambda^{0} \subset \Lambda$ and an 
uncountable dense subset $\Cal G \subset \Delta$ such that: 
\roster 
\item $\Cal S := \Delta - \Cal G$ is countable and dense; 
\item $\Lambda^{0} \subset \Lambda_{t}$ for all $t \in \Delta$;
\item $\Lambda^{0} = \Lambda_{t}$ for all $t \in \Cal G$; 
\item $\Lambda^{0} \not= \Lambda_{t}$ for all $t \in \Cal S$. 
\endroster 
We call a point $t \in \Cal G$ {\it generic} and a point $s \in \Cal S$ 
{\it special}. Special points are nothing but the points $s$ at which 
the Picard numbers $\rho(\Cal X_{s})$ jump above. 
\par
\vskip 4pt
Under this notation, we can state our main result as follows:
\proclaim{Main Theorem} Assume that $\varphi : \Cal X \rightarrow \Delta$ 
is projective, i.e. there exists a $\varphi$-ample invertible sheaf $\Cal L$. 
Then:
\roster 
\item There exist a (possibly empty) finite subset $\Cal F \subset \Cal S$, a 
group $G^{0}$, and a positive integer $N$ depending only on $\varphi$ such that 
$$G^{0} < \text{Aut}(\Cal X_{t})\,\, \text{for all}\,\, t \in \Delta - \Cal F$$ 
and 
$$[\text{Aut}(\Cal X_{t}) : G^{0}] \leq N \,\, \text{for all}\,\, t \in \Cal G.$$ 
In particular, the map
$$\text{Aut} : \Delta \rightarrow \{\text{groups}\}/\equiv\, ;\, 
t \mapsto [\text{Aut}(\Cal X_{t})]$$ 
is "upper-semicontinuous" with respect to the co-finite topology of 
$\Delta -\Cal F$. 
\item There exists a subset $\Cal D \subset \Cal S$ such that 
$\Cal D$ is dense in $\Delta$ and that
$$\vert \text{Aut}(\Cal X_{t}) \vert = \infty\,\, \text{for all}\,\, t \in 
\Cal D.$$ 
\item There exists a projective family of K3 surfaces $\varphi : \Cal X 
\rightarrow \Delta$ such that $\Cal F \not= \emptyset$.  
\endroster 
\endproclaim 
As a direct consequence of the main Theorem, one obtains the following:
\proclaim{Corollary} Let $f : \Cal X \rightarrow \Delta$ be a (not necessarily projective) non-trivial family of K3 surfaces. Then, there is a dense subset $\Cal D \subset \Delta$ such that
$\vert \text{Aut}(\Cal X_{t}) \vert = \infty$ for all $t \in \Cal D$. In 
particular, the nef cone $\overline{A}(\Cal X_{t})$ is not finite rational 
polyhedral if $t \in \Cal D$.
\endproclaim 
Compare the results with certain boundedness of K3 surfaces with $\vert 
\text{Aut}(X) \vert < \infty$ ([Ni6]). We shall prove the main Theorem and Corollary 
in Section 1.
\par
\vskip 4pt
The assertion (2) and Corollary are somewhat surprizing. For instance, let us 
take the component $\Cal H$ of the Hilbertscheme consisting of quartic 
K3 surfaces and consider the universal family $u : \Cal U \rightarrow \Cal H$. 
Then 
for any sufficiently general $\Delta \rightarrow \Cal H$, the induced family $\varphi : \Cal X \rightarrow \Delta$ satisfies 
$\text{Pic}(\Cal X_{t}) = {\bold Z}\Cal L_{t}$ for $t \in \Cal G$, 
where $\Cal L_{t}$ is the plane section class of $\Cal X_{t} \subset 
\bold P^{3}$. So $\vert \text{Aut}(\Cal X_{t}) \vert < \infty$ if $t \in 
\Cal G$. A bit more strongly, there also exists $M$ with $\vert 
\text{Aut}(\Cal X_{t}) \vert < M$ for all $t \in \Cal G$ by the 
boundedness of the finite subgroups of automorphism groups of K3 surfaces 
([Ni1], [Mu] see also [Kn2]). The statement (2) claims, however, 
that there exists a dense subset $\Cal D$ such that 
$\vert \text{Aut}(\Cal X_{t}) \vert = \infty$ for all 
$t \in \Cal D$ even in this family. This also provides an explicit example of a family in which the automorphism groups actually jump above. On the other hand, one has by Koll\'ar([Bo]): {\it If $f : \Cal Y \rightarrow \Delta$ is a family 
of Calabi-Yau manifolds in $\vert -K_{V} \vert$ of a Fano manifold $V$, then 
the nef cones $\overline{A}(\Cal Y_{t})$ are finite rational 
polyhedral, whence $\vert \text{Aut}(\Cal Y_{t}) \vert < \infty$, for all 
$t \in \Delta$ provided that $\dim V \geq 4$}. Thus, the statements similar to Theorem (2) and Corollary do not hold for Calabi-Yau manifolds of higher 
dimension. 
\par
\vskip 4pt 
The assertion (1) mostly justifies the phrase quoted at the beginning, 
while the assertion (3) denies the phrase in the most strict sense. By the statement (1) and Corollary, one might expect something more geometric 
behind such as a 
covariant relation among jumping of automorphism groups, jumping of Picard 
numbers, those of smooth rational curves and those of elliptic pencils 
in a family of K3 surfaces. However, there are no such relations in general. 
Indeed, the second aim of this note is to point out the following: 

\proclaim{Proposition} There exist families of projective K3 surfaces 
$\varphi^{i} : \Cal X^{i} \rightarrow \Delta$ ($1 \leq i \leq 5$) such that the central fiber $\Cal X_{0}^{i}$ satisfies that $\rho(\Cal X_{0}^{i}) = 19$, 
$\vert \text{Aut}(\Cal X_{0}^{i})\vert < \infty$, contains finitely 
many smooth rational curves and admits finitely many elliptic pencils, but 
that generic fiber $\Cal X_{t}^{i}$ satisfies $\rho(\Cal X_{t}^{i}) = 3$ and 
that:
\roster 
\item $\vert \text{Aut}(\Cal X_{t}^{1})\vert = \infty$ and $\Cal X_{t}^{1}$ 
contains infinitely many smooth rational curves but admits no elliptic pencil. 
\item $\vert \text{Aut}(\Cal X_{t}^{2})\vert = \infty$ and $\Cal X_{t}^{2}$ 
contains no smooth rational curve but admits infintely many elliptic pencils. 
\item $\vert \text{Aut}(\Cal X_{t}^{3})\vert = \infty$ and $\Cal X_{t}^{3}$ 
has infinitely many smooth rational curves and infinitely many elliptic 
pencils.  
\item $\vert \text{Aut}(\Cal X_{t}^{4})\vert = \infty$ and $\Cal X_{t}^{4}$ 
has neither smooth rational curve nor elliptic pencil. 
\item $\vert \text{Aut}(\Cal X_{t}^{5})\vert < \infty$ and $\Cal X_{t}^{5}$ 
has both (necessarily finitely many) smooth rational curves and elliptic 
pencils. 
\endroster 
\endproclaim 
We shall construct these families in Section 2 by deforming a very remarkable 
K3 surface found by Nikulin [Ni3] (See also [Kn1]). In our construction, the 
morphisms $\varphi^{i}$ are not projective. 
\head
{Acknowledgement}
\endhead
Main part of work has been done during the author's stay at Bayreuth in 
March 2001 under financial support by Germany. Idea of work has been brought 
to the author by Professor Thomas Peternell through their discussion on 
another problem (in progress). The author would like to express his best thanks to him for his invitation to Bayreuth as well as valuable discussion. 

\head
{1. Proof of the main Theorem} 
\endhead 
\demo{Proof of Theorem (1)} Set 
$$l := \tau_{t}([\Cal L_{\Cal X_{t}}]) \in \Lambda_{t}\, (\, \subset 
\Lambda\,).$$ 
Since $l$ is independent of $t$, we have $l \in \Lambda^{0}$. We may assume 
that the class $l$ is primitive. Let $A(\Cal X_{t})$ be the ample cone of $\Cal X_{t}$. We set: 
$$A_{t} := \tau_{t}(A(\Cal X_{t}))\, \subset \Lambda_{t} \otimes \bold R;$$ 
$$A_{t}^{0} := A_{t} \cap (\Lambda^{0} \otimes \bold R);$$
$$G_{t} := \{g \in O(\Lambda^{0}) \vert\,\, g \vert(\Lambda^{0*}/\Lambda^{0}) = id\, , \, g(l) \subset A_{t}^{0}\},$$ 
for each $t \in \Delta$. Here $\Lambda^{0*}$ is the dual lattice of 
$\Lambda^{0}$ and is regarded as an overlattice of $\Lambda^{0}$ 
via the cup product $(*, **)$, which is non-degenerate on $\Lambda^{0}$ by the 
Hodge index Theorem. Note that $l \in A_{t}^{0}$ for each $t$. For simplicity of 
description, we often identify $\Lambda_{t}$ with $\text{NS}(\Cal X_{t})$ by $\tau_{t}$. For instance, we say that an element $a \in A_{t}$ is ample on 
$\Cal X_{t}$. The next Claim completes the proof of Theorem (1): \enddemo
\proclaim{Claim 1} Choose $p \in \Cal G$. Then, there exists a (possibly empty) finite subset $\Cal F \subset \Cal S$ such that
\roster 
\item $A_{t}^{0} \subset A_{p}$ for all $t \in \Delta$, and $A_{t} = A_{p}$ for 
all $t \in \Cal G$. 
\item $G_{t} \supset G_{p}$ for all $t \in \Delta - \Cal F$, and $G_{t} = G_{p}$ for all $t \in \Cal G$. 
\item There is a natural embedding $\iota_{t} : G_{t} \hookrightarrow \text{Aut}(\Cal X_{t})$ for all $t \in \Delta$. Moreover, there exists a positive integer $N$ such that $[\text{Aut}(\Cal X_{t}) : \iota_{t}(G_{t})] \leq N$ for all $t 
\in \Cal G$. 
\endroster
\endproclaim
\demo{Proof of Claim 1} The first assertion is clear if $\text{rank}\, \Lambda^{0} = 1$. 
We assume $\text{rank}\, \Lambda^{0} \geq 2$. Let $h \in A_{t}^{0}$. Then $(h^{2}) > 0$ and $(h,l) > 0$. If $h$ is not ample on $\Cal X_{p}$, then there is an 
irreducible curve $C$ on $\Cal X_{p}$ such that $([C], h) \leq 0$. This is due 
to the real version of the Nakai-Moishezon criterion of ampleness by Campana 
and Peternell [CP]. This $C$ must be a smooth rational curve by the Hodge 
index Theorem and the adjunction formula. Moreover, since $C$ is a curve on 
$\Cal X_{p}$, we have $[C] \in \Lambda^{0}$. Thus $[C] \in \Lambda_{t}$ as 
well. Since $(C^{2}) = -2$, either $-[C]$ or $[C]$ is represented by a 
(non-zero) effective curve on $\Cal X_{t}$ by the Riemann-Roch Theorem. Since 
$([C], l) > 0$ by the 
ampleness of $l$ plus effectiveness of $C$ on $\Cal X_{p}$ and since $l$ is also ample on $\Cal X_{t}$, it is $[C]$ that is represented by an effective curve 
on $\Cal X_{t}$. ({\it Here the assumption $\varphi$ being projective is 
essential.}) However, since $h \in A_{t}$, one would then have $([C], h) > 0$, 
a contradiction to the previous inequality. {\it Here we used the fact that the 
values} $(*, **)$ {\it are independent on which K3 surfaces} $\Cal X_{t}$ {\it we are 
arguing.} Therefore we have $A_{t}^{0} \subset A_{p}$ for all $t \in \Delta$. ({\it Note, however, that the other inclusion is false in general as we will see 
in the proof of (3).}) Since $A_{t}^{0} = A_{t}$ if $t \in \Cal G$, replacing $p$ by $t$, one has 
$A_{t} = A_{p}$ for all $t \in \Cal G$. Now, by the definition of $G_{*}$, 
one has $G_{t} = G_{p}$ for all $t \in \Cal G$ as well. 
\par 
\vskip 4pt 
We shall show the first assertion of (2). From now, we write $A_{p}$ by $A$. 
Set 
$$\Cal F := \{t \in \Cal S \vert G_{p} \not\subset G_{t}\}.$$  
We need to show that $\Cal F$ is finite. Recall that the group $G_{p}$ is 
finitely generated by [BC] (See also [St]). Indeed, $G_{p}$ is a quotient of an arithmetic subgroup of $O(\Lambda^{0})$. Choose a set of generators of $G_{p}$ and denote it by  $\{g_{i} \vert 1 \leq i \leq n \}$. Set for each $i$: 
$$\Cal F_{i} := \{ t \in \Cal S \vert g_{i} \not\in G_{t}\}.$$ 
It suffices to show that $\vert \Cal F_{i} \vert < \infty$. Choose $i$ and put $g := g_{i}$. Assume that $t \in \Cal F_{i}$. 
Then, by the definition, $g(l) \not\in A_{t}$. (Note that $g(l) \in 
\Lambda^{0}$ by $g \in G_{p}$.) Since $(g(l)^{2}) = (l^{2}) > 0$ and 
$(g(l), l) > 0$ by $l$, $g(l) \in A$, we see that $g(l)$ is in the positive 
cone of $\Cal X_{t}$ by the Hodge index Theorem. Then, 
applying the Riemann-Roch Theorem, one finds an effective divisor $D$ on 
$\Cal X_{t}$ such that $g(l) = [D]$. Since $g(l)$ is not ample on a K3 surface 
$\Cal X_{t}$, there exists a smooth rational curve $C_{t}$ such that 
$(g(l), C_{t}) \leq 0$. 
\par 
\vskip 4pt
Let $\tilde{g}$ be the element of $O(\Lambda)$ such that 
$\tilde{g} \vert \Lambda^{0} = g$ and $\tilde{g} \vert \Lambda^{0 \perp} 
= id$, where $\Lambda^{0\perp}$ is the orthogonal complement of $\Lambda^{0}$ 
in $\Lambda$. Such an extension exists, because $g \vert \Lambda^{0 *}/\Lambda^{0} = id$. Since 
$T_{\Cal X_{t}} \subset \Lambda^{0 \perp}$, where $T_{\Cal X_{t}}$ is the 
transcendental lattice of $\Cal X_{t}$, this $\tilde{g}$ is also a Hodge 
isometry of $\Cal X_{t}$. Thus, $E := \tilde{g}^{-1}([C_{t}])$ is an 
integral $(1,1)$-class on $\Cal X_{t}$ and satisfies $([E]^{2}) = -2$. In 
particular, $[E]$ or $-[E]$ is represented by an effective curve on 
$\Cal X_{t}$ by the Riemann-Roch Theorem. Therefore, $(g(l), C_{t}) = (l, E) 
\not= 0$ by the ampleness of $l$, and we have $(g(l), C_{t}) < 0$. Then $C_{t} \subset \text{Supp}\,D$ and one has $(l, C_{t}) \leq (l, g(l))$. Note that $(l, g(l))$ is independent on $t$. So, if there are 
infinitely many $t$ such that $t \in \Cal F_{i}$, then by [Gr], there exists a 
component $\Cal H$ of the relative Hilbertscheme 
$\text{Hilb}_{\Cal X/\Delta}^{l}$ to which infinitely many of $C_{t}$ 
belong. (See also a very nice book [Kl] for the Hilbertschemes.) 
This $\Cal H$ must dominate $\Delta$ and one has then an 
effective curve $C$ on 
$\Cal X_{p}$ such that  $([C], g(l)) = ([C_{t}], g(l)) < 0$ (for $C$ and 
$C_{t}$ belong to the same component $\Cal H$). However, 
this contradicts the ampleness of $g(l)$ on $\Cal X_{p}$. Hence 
$\vert \Cal F_{i} 
\vert < \infty$ for each $i$ and we are done for (2).  
\par 
\vskip 4pt 
We shall show the assertion (3).  As before, each element $g_{t} \in G_{t}$ 
extends a Hodge isometry $\tilde{g}_{t}$ of $\Cal X_{t}$ such 
that $\tilde{g}_{t} \vert \Lambda^{0} = g_{t}$ and $\tilde{g}_{t} \vert \Lambda^{0 \perp} = id$. This $\tilde{g}_{t}$ is also effective, because 
$$\tilde{g}_{t}(l)\, = \, g_{t}(l) \, \in \, A_{t}^{0} \, \subset \, A_{t}.$$ 
There then exists a unique $f_{t} \in 
\text{Aut}(\Cal X_{t})$ such that $f_{t}^{*} \vert \Lambda^{0} = g_{t}$ and 
$f_{t}^{*} \vert \Lambda^{0 \perp} = id$ by the global Torelli Theorem for 
K3 surfaces. The map $\iota_{t}$ is given by $g_{t} \mapsto f_{t}$. Now the 
assertion follows from the next facts:
\roster
\item $\tau_{t}(\text{NS}(\Cal X_{t})) = \Lambda^{0}$ for all $t \in \Cal G$; 
\item $\Lambda^{0 *}/\Lambda^{0}$ is finite (and is independant of $t$);
\item The subgroup $\text{Ker}(\text{Aut}(\Cal X_{t}) \rightarrow O(NS(\Cal X_{t})^{\perp}))$ of $\text{Aut}(\Cal X_{t})$ 
is of finite index less than or equal to $66 \times \vert \text{Aut}(\Lambda^{0*}/\Lambda^{0}) \vert$ if $t \in \Cal G$ ([Ni1], see also [St]). \qed
\endroster 
\enddemo
Now we are done for Theorem (1). \qed 
\demo{Proof of Theorem (2)} We may consider the case 
$\vert \text{Aut}(\Cal X_{t}) \vert < \infty$ for all $t \in \Cal G$. (If otherwise, one has  $\vert \text{Aut}(\Cal X_{t}) \vert = \infty$ for all 
$t \in \Delta - \Cal F$ by Theorem (1) and may put 
$\Cal D = \Cal S - \Cal F$.) Since $\Cal G$ is dense, the statement will 
follow from the next: \enddemo
\proclaim{Claim 2} For any given $p \in \Cal G$ and for any open neighbourhood 
$p \in U$, there exists $t \in \Cal S \cap U$ such that 
$\vert \text{Aut}(\Cal X_{t}) \vert = \infty$. 
\endproclaim 
\demo{Proof of Claim 2} By abuse of notation, we shall write again 
$U = \Delta$ and $\varphi \vert U = \varphi$. Assuming $\vert 
\text{Aut}(\Cal X_{t}) \vert < \infty$ for all 
$t \in \Delta$, 
we shall derive a contradiction. Set $r_{0} := \text{rank}\, \Lambda^{0}$. 
By shrinking, we may assume that $\Cal F = \emptyset$. We shall argue 
dividing into the following two cases (1) and (2) (by shrinking again, we see that these 
two cover all the possible cases): 
\roster
\item There exists a sequence $\{t_{n}\}_{n \geq 1}\, (\subset \Cal S)$ such 
that $\text{rank}\, \Lambda_{t_{n}} \geq 3$ for all $n$ and that $\lim_{n \rightarrow \infty} 
t_{n} = p$;
\item $r_{0} = 1$ and $\text{rank}\, \Lambda_{t} = 2$ for all $t \in \Cal S$.
\endroster 
\demo{Case (1)} Since $\text{rank}\, \Lambda_{t_{n}} \geq 3$ and $\vert 
\text{Aut}(\Cal X_{t_{n}}) \vert < \infty$, by [PSS], $\Cal X_{t_{n}}$ 
contains 
only 
finitely many smooth rational curves, say, 
$$C_{n1}\, ,\, C_{n2}\, ,\, \cdots \, , \, C_{nk(n)}.$$ 
By Kovacs [Kv], the classes $[C_{ni}]$ generate the Kleiman-Mori cone 
$N_{n} (:= 
\overline{NE}(\Cal X_{t_{n}}))$ of 
$\Cal X_{t_{n}}$. Recall also that the isomorphism classes of both the 
Picard lattices and the dual graphs of smooth rational curves, or in other 
words, the fundamental domains of the $2$-reflection groups, of K3 surfaces with $\vert \text{Aut}(X) \vert < \infty$ and $\rho(X) \geq 3$ are finite. This deep result is due to Nikulin [Ni 3, 4, 5] and it is this finiteness that we need the case assumption essentially. ({\it Note that such finiteness is clearly false if $\rho \leq 2$.}) Therefore we may assume that $\rho := \text{rank} 
\Lambda_{t_{n}}$ and $k := k(t_{n})$ are constant and the intersection 
matrices $(C_{ni}, C_{nj})$ are also independent of $n$. ({\it However, we do 
not know 
a priori whether the classes $[C_{ni}]$ form constant subsystems of $\Lambda \times \Delta$.}) 
\par 
\vskip 4pt 
Consider first the case of $r_{0} \geq 3$. Then $\Cal X_{p}$ has finitely many 
smooth rational curves as well by $\vert \text{Aut}(\Cal X_{p}) \vert < 
\infty$. Denote all of them by $D_{j}$ ($1 \leq j \leq k'$). Then one has $l = \sum_{j = 1}^{k'} a_{j}D_{j}$, 
where $a_{j}$ are non-negative rational numbers. Since $[D_{j}] \in \Lambda^{0} \subset \Lambda_{t}$, each $D_{j}$ extends $\bold P^{1}$-bundles $\Cal D_{j} 
\rightarrow \Delta$ near $p$. (See for instance [Og].) Thus, for large $n$, one has $\Cal D_{jt_{n}} 
\in \{C_{ni}\}_{i = 1}^{k}$ and $l =  \sum_{j = 1}^{k'} a_{j}\Cal D_{jt_{n}}$. 
Then the set $\{(C_{ni}, l)\}_{n \geq 1}$ is bounded for all $i$ by the 
finiteness of dual graphs. Thus, using the existence of relative 
Hilbertscheme and arguing as before, one finds effective divisors $C_{0i}$ ($i = 1, 2, \cdots , k$) on $\Cal X_{p}$ whose dual graph is same as $C_{ni}$. However, one would then have $\text{rank}\, 
\Lambda^{0} \geq \text{rank}\, \Lambda_{t} > r_{0}$, a contradiction. 
\par 
\vskip 4pt 
Next consider the case of $r_{0} = 2$. By $\vert \text{Aut}(\Cal X_{p}) \vert 
< \infty$, the Kleiman-Mori cone $N_{p}$ of $\Cal X_{p}$ is generated by 
either (1) two 
classes of smooth rational curves, (2) two classes of elliptic pencils, (3) 
the class of a smooth rational curve and the class of elliptic pencil, or (4) 
the class of a smooth rational curve and an irrational ray $\bold R_{\geq 0} e$ with $(e^{2}) = 0$ (See for instance [Kv]). In the first three cases, the same argument as above immediately gives a contradiction if we notice the following facts:
\roster 
\item Each elliptic pencil on $\Cal X_{p}$ extends a family of elliptic pencils near $p$. This is because $\Lambda^{0} \subset \Lambda_{t}$. (See 
for instance [Og]).
\item Each elliptic pencil on $\Cal X_{t_{n}}$ has at least one fiber 
consisting of smooth rational curves. This follows from Shioda's formula 
of the Mordell-Weil rank of the Jacobian (see 
for instance [Sh]) and 
the assumption $\text{rank}\, \Lambda_{t_{n}} \geq 3$. Indeed, the 
Mordell-Weil rank must be $0$ by $\vert \text{Aut}(\Cal X_{t_{n}}) \vert < 
\infty$. 
\endroster 
So, we may consider the last case (4). Since the smooth rational curve $C$ on 
$\Cal X_{p}$ extends sideways (by $[C] \in \Lambda^{0}$), we may assume that 
$[C] = [C_{n1}]$ for large $n$. This $[C]$ generates one of the two boundary 
rays 
of $N_{p}$. 
Since $A_{t_{n}}$ is finite rational polyhedral, so is $A_{n}^{0} := 
A_{t_{n}} \cap (\Lambda^{0} \otimes \bold R)$. Therefore, by the case 
assumption, one has $A_{n}^{0} \not= A_{p}$ (while one has now $A_{n}^{0} \subset A_{p}$ 
by Claim 1 (1).) There then exist a nef and big divisor $H_{n}$ on $\Cal X_{t_{n}}$ and a smooth rational curve $C_{n2}$ on $\Cal X_{t_{n}}$ such that 
$(H_{n}, C_{n2}) 
= 0$ and that $[H_{n}] \in A_{p}$. (Note that $([H_{n}], [C_{n1}]) = 
([H_{n}], [C]) > 0$ because $[H_{n}] \in A_{p}$.) Choose a euclidean norm $\Vert * \Vert$ on $\Lambda \otimes \bold R$. If the set $\{\Vert C_{n2} \Vert\}_{n \geq 1}$ is bounded, 
then, since there are only finitely many lattice points with bounded norm, 
the set of intersection numbers $\{(C_{n2}, l)\}_{n \geq 1}$ is bounded as 
well. 
Then, by the existence of the relative Hilbertscheme, $\Cal X_{p}$ contains 
an effective divisor $D_{02}$ such that $[D_{02}] = [C_{n2}]$. However, the class 
$[D_{02}]$ gives then another rational boundary of $N_{p}$, a 
contradiction to the case assumption. Therefore, 
$\{\Vert C_{n2} \Vert\}_{n \geq 1}$ is unbounded. Then, passing to a 
subsequence, we have $\lim_{n \rightarrow \infty} \Vert C_{n2} \Vert = 
\infty$. Set $x_{n2} := C_{n2}/\Vert C_{n2} \Vert$. Since 
$\Vert x_{n2} \Vert = 1$, by passing to some subsequence, we find the limit 
$x_{02} := \lim_{n \rightarrow \infty} x_{n2}$ in $\Lambda \otimes 
\bold R$. We calculate that $\Vert x_{02} \Vert = 1$ and that 
$$(x_{02}^{2}) = \lim_{n \rightarrow \infty} -2/\Vert C_{n} \Vert = 0.$$ 
Moreover, we calculate that 
$$(x_{02}, [\omega_{\Cal X_{p}}]) = \lim_{n \rightarrow \infty}
(x_{n2}, [\omega_{\Cal X_{t_{n}}}]) 
= 0.$$ 
Here $\omega_{*}$ is the $2$-form on $\Cal X_{*}$ given by the base change of 
a local section $\omega$ of $\varphi_{*}\Omega_{\Cal X/\Delta}^{2}$ around $p$. 
Thus, $x_{02}$ is also a real $(1,1)$-class on $\Cal X_{p}$. Set $h_{n} := H_{n}/\Vert H_{n} 
\Vert$. Since $h_{n} \in A_{p}$, passing to a subsequence, we have the limit
$$h_{0} := \lim_{n \rightarrow \infty} h_{n} \in \overline{A}_{p}.$$ 
Since $(h_{n},x_{n2}) = 0$, one has $(h_{0}, x_{02}) = 0$ as well. Then, by 
the Hodge index Theorem, one has that $x_{02} = \alpha h_{0}$ and $x_{02} \in 
\overline{A}_{p}$. In particular, $\bold R_{>0}x_{02}$ gives another boundary of $N_{p}$. 
On the other hand, since $(C_{n1}, C_{n2})$ is constant, one calculates 
$$([C], x_{02}) = \lim_{n \rightarrow \infty}([C_{n1}], [C_{n2}])/\Vert C_{n2} 
\Vert 
= 0.$$ 
Since $\Lambda^{0}$ is of rank $2$, this implies $\bold R x_{02} = 
[C]^{\perp}$. However, $\bold R_{>0} x_{02}$ would be then a rational boundary, a contradiction to the case assumption. 
\par
\vskip 4pt 
Finally consider the case where $r_{0} = 1$. If 
$\{\Vert C_{ni} \Vert\}_{n \geq 1}$ is bounded for some $i$, then, as before, 
$\{(C_{ni}, l)\}_{n \geq 1}$ is bounded and we have an effective 
curve 
$D$ such that $[D] = [C_{ni}]$ for large $n$. This class $[D]$ satisfies 
$[D] \in \text{NS}(\Cal X_{p})$ and $([D]^{2}) = -2$. However, this is impossible, because $\Cal X_{p}$ is projective and $\text{rank}\, \text{NS}(\Cal X_{p}) 
= r_{0} = 1$. Therefore the sets $\{\Vert C_{ni} \Vert\}_{n \geq 1}$ are 
unbounded for each $i$. 
We set $x_{ni} := C_{ni}/\Vert C_{ni} \Vert$. Passing to a subsequence, we 
have 
$\lim_{n \rightarrow \infty} \Vert C_{ni} \Vert = \infty$ and get the limit 
points $x_{0i} := \lim_{n \rightarrow \infty} x_{ni}$. Since $(C_{ni}, C_{nj})$ are all bounded, we have that $(x_{0i}^{2}) = (x_{0i}, x_{0j}) = 0$ for all $i$ and $j$. Moreover, by the same argument as in the previous case, we see that 
$x_{0i}$ are all real $(1,1)$-classes on $\Cal X_{p}$. Then, by the Hodge 
index Theorem, there exists $e$ with $(e^{2}) = 0$ such that $x_{0i} \in 
\bold R e$ for all $i$. Then, the Kleiman-Mori cones $N_{n}$, which are finite 
polyhedral cones generated by $x_{ni}$ ($1 \leq i \leq k$), accumulate to the 
ray $\bold R e$ 
if $n \rightarrow \infty$. However, since $l \in N_{n}$, we would then have 
$l \in \bold R e$, a contradiction to $(l^{2}) > 0$. Now we are done in Case 1. \qed \enddemo
\demo{Case (2)} Observe that if $\text{rank}\, \Lambda_{t} = 2$, then 
$\vert \text{Aut}(\Cal X_{t}) \vert = \infty$ if and 
only if the quadratic form of $\Lambda_{t}$ represents neither $0$ nor $-2$. 
Let us consider the period map
$$\pi : \Delta \rightarrow \Cal P := \{[\omega] \in \bold P(\Lambda \otimes 
\bold C) \vert (\omega, \omega) = 0, (\omega, \overline{\omega}) > 0\}.$$ 
For each $a \in \Lambda$, we set 
$$H_{a} := \{[\omega] \in \Cal P \vert (a, \omega) = 0\} \subset \Cal P.$$ 
Then $\pi(\Delta) \subset H_{l}$. Moreover, by the case assumption, we have that $t \in \Cal S$ if and only if there 
exists $a \in \Lambda - \bold Z l$ such that $\pi(t) \in H_{a}$ and that 
$\Lambda_{t} = \bold Z \langle l, a \rangle$. (Recall that $l$ is assumed to 
be primitive.) For this $t$, one has the following: \enddemo
\proclaim{Claim 3} There exist an element $h \in \Lambda$ and arbitrarily 
large positive integers $n$ such that $\bold Z \langle l, na + h \rangle$ is 
primitive in $\Lambda$ and represents neither $0$ nor $-2$. \endproclaim 
\demo{Proof of Claim 3} Set $(l, l) = 2A$, $(l, a) = B$ and $(a, a) = 2C$. Since the 
primitive embedding $\Lambda_{t} \hookrightarrow \Lambda = U^{\oplus 3} \oplus 
E_{8}(-1)^{\oplus 2}$ is unique up to isomorphism by the primitive embedding 
Theorem [Ni2, Theorem 1.14.4], we may assume that 
$$l = e_{11} + Ae_{12}\,\, , \,\, a = Be_{12} + e_{21} + Ce_{22}\, ,$$ 
where $e_{i1}$ and 
$e_{i2}$ ($1 \leq i \leq 3$) are the standard basis of $U^{\oplus 3}$. Set 
$h := e_{31} - me_{32}$. Then $(h^{2}) = -2m$ and $(h, l) = (a, l) = 0$. 
Put:
$$L_{m,n} := \bold Z \langle l, na + h \rangle.$$ 
Then this lattice $L_{m,n}$ is primitive in 
$\Lambda$ and its quadratic form is 
$$Q(x, y) := 2\{Ax^{2} + nBxy + (n^{2}C - m)y^{2}\}.$$ 
If $A \geq 2$, we take $n = AN$ and $m = AM$. Then $Q(x, y)$ is divisible by $2A$. 
Thus $L_{m,n}$ does not represent $-2$. Moreover, the discriminant of $Q$ is 
$$4A^{2}\{N^{2}(B^{2} - 4AC) + 4M\}.$$ 
This is not square if $M$ and $N$ are general. Thus $L_{m,n}$ does not 
represent $0$, either. If $A = 1$, we set $n = 4N$ and $m = 8M$. Then 
$Q(x, y)$ and its discriminant are 
$$Q(x, y) = 2x^{2} + 8NBxy + 8(4N^{2}C - M)y^{2};$$
$$D := 64\{N^{2}(B^{2} - 4C) + M\}.$$ 
Then, $Q(x, y) \equiv 2$ or $0$ modulo $8$. Thus $L_{m, n}$ does not 
represent $-2$. Moreover, $D$ is not sqaure if $M$ and $N$ are general. 
Thus $L_{m,n}$ does not represent $0$, either. \qed \enddemo
 
Let us return back to the original situation. If the lattice 
$\bold Z \langle l, a \rangle$ (for $t$) represents neither $0$ nor $-2$ then 
we are 
done. If otherwise, using Claim 3, we find an element $h \in \Lambda$ and 
an element $n \in \bold Z$ such that 
$\bold Z \langle l, na + h \rangle$ is primitive in $\Lambda$ and represents 
neither $0$ nor $-2$. Since $H_{na + h} = H_{a + (h/n)}$, the hyperplane 
$H_{na + h}$ 
meets $\pi(\Delta)$ at $\pi(t')$, where $t'$ is a point very close to $t$, 
if $n$ is 
sufficiently large. For this $t'$, one has $\Lambda_{t'} = 
\bold Z \langle l, na + h \rangle$ and $\vert \text{Aut}(\Cal X_{t'}) \vert 
= \infty$. \qed 
\par
\vskip 4pt
Now we are done for Theorem (2). \qed 
\demo{Proof of Theorem (3)} Let $T$ be a K3 surface such that 
$\text{NS}(T) = U \oplus A_{1}(-1)$ and $\vert \text{Aut}(T) \vert < \infty$. 
Such a K3 surface exists by [Ni3]. 
Let us consider the Kuranishi family $u : \Cal U \rightarrow \Cal K$ of $T$.
Choosing a marking 
$$\tau : R^{2}u_{*}\bold Z_{\Cal U} \rightarrow \Lambda \times \Cal K,$$ 
one can define the period map 
$$\pi : \Cal K \rightarrow \Cal P.$$ 
Since $\pi$ is a local isomorphism by the local Torelli Theorem for K3 
surfaces, we may identify $\Cal K$ with a small open set of the
period domain $\Cal P$ (denoted again by $\Cal P$) through $\pi$. Let $\Pi$ be 
a general primitive sublattice of $\text{rank}\,2$ of $U \oplus A_{1}(-1)$. Then $\Pi$ represents neither 
$0$ nor $-2$. Let us take the subset $\Cal K^{0} \subset \Cal K = \Cal P$ 
defined by the equations 
$$(a, \omega) = (b, \omega) = 0,$$ 
where $a$ and $b$ form an integral basis of $\Pi$. Consider the induced 
family $u^{0} : \Cal U^{0} \rightarrow \Cal K^{0}$. Here $\Cal K^{0}$ is of 
dimension $18$. 
By construction, we have $0 \in \Cal K^{0}$ and $\text{NS}(\Cal U_{t}^{0}) \simeq \Pi$, whence $\vert \text{Aut}(\Cal U_{t}^{0}) \vert = \infty$ for generic 
$t$. Then, 
the induced family $\varphi : \Cal X \rightarrow \Delta$ for a generic linear 
section $0 \in 
\Delta \subset \Cal K^{0}$ satisfies $0 \in \Cal F$. In this example, we have 
$A_{0}^{0} \not= A_{t}$ for $t$ being generic. Indeed, $\partial A_{0}^{0} = 
\partial (A_{0} \cap \Pi)$ is rational, while $\partial A_{t}$ is irrational.  
\qed \enddemo
Now we are done. Q.E.D. of the main Theorem. \qed \enddemo 
\demo{Proof of Corollary}  
Since $\Lambda^{0} \subset \Lambda_{t}$, an ample divisor on $\Cal X_{p}$ 
extends 
sideways (as a line bundle) if $p \in \Cal G$. In particular, $f$ is locally projective around $p$. 
The result now follows from the Theorem together with the density of $\Cal G$. 
\qed \enddemo
\head
{2. Proof of Proposition} 
\endhead   
We shall construct explicit examples with required properties by deforming a 
remarkable K3 surface found by Nikulin and Kondo ([Ni3], [Ko1]). This is a K3 
surface $S$ which satisfies the following properties:
\roster 
\item $\rho(S) = 19$ and $\text{NS}(S) = U \oplus E_{8}(-1) \oplus E_{8}(-1) 
\oplus A_{1}(-1)$ as a lattice, 
\item $S$ contains exactly $24$ smooth rational curves and finitely many 
($> 0$) elliptic pencils; 
\item $\text{Aut}(S)$ is a finite group (and is indeed isomorphic to 
$S_{3} \times 
\mu_{2}$). 
\endroster 
As before, let us take the Kuranishi family $u : \Cal U \rightarrow \Cal K$ of 
$S$ and choose a marking 
$$\tau : R^{2}u_{*}\bold Z_{\Cal U} \rightarrow \Lambda \times \Cal K.$$ 
Then $S = \Cal U_{0}$ and one has the period map 
$$\pi : \Cal K \rightarrow \Cal P.$$
Again as before, we identify $\Cal K$ with a small open set of the
period domain $\Cal P$ (denoted again by $\Cal P$) through $\pi$. Let 
$$e\, , \, f\, , \, v_{ij}\, (i = 1\, ,\, 2\, , \, 1 \leq j \leq 8)\, , \, v$$ 
be the elements of $\Lambda$ corresponding to the standard integral basis of 
$\text{Pic}(S) = U \oplus E_{8}(-1) \oplus E_{8}(-1) \oplus A_{1}(-1)$. 

\demo{Proof of (1)} Let us take the subset $0 \in \Cal K^{1} \subset \Cal K = 
\Cal P$ defined by the equations 
$$(e + 3nf, \omega) = (v_{11}, \omega) = (v_{21}, \omega) = 0,$$ 
where $n$ is an integer non-divisible by $3$. Consider the induced family 
$u^{1} : \Cal U^{1} \rightarrow \Cal K^{1}$. Here $\Cal K^{1}$ is of 
dimension 
$17$. By construction, $\text{NS}(\Cal U_{t}^{1}) \simeq \langle e + 3nf, 
v_{11}, v_{21} \rangle$ and the intersection matrix is: 
$$\pmatrix 6n & 0 & 0\\
            0 & -2 & 0\\
            0 & 0 & -2\\ 
\endpmatrix ,$$
if $t \in \Cal K^{1}$ is generic. The discriminant of this matrix is $24n$. 
Since the quadratic form 
$$Q_{1}(x, y) := 6nx^{2} - 2y^{2} - 2z^{2}$$ 
does not represent $0$, $\Cal U_{t}^{1}$ has no elliptic pencil. Since there are only finitely many 
isomorphism classes of lattices of rank $3$ which are isomorphic to the Picard 
lattices of K3 surfaces $X$ with $\vert \text{Aut}(X) \vert < \infty$ [Ni4], 
one has $\vert \text{Aut}(\Cal U_{t}^{1}) \vert = \infty$ for $n$ being 
sufficiently large. Moreover, $\Cal U_{t}^{1}$ admits at least one smooth 
rational curve. This follows from the fact that $Q_{1}(x, y)$ represents $-2$ 
plus the Riemann-Roch Theorem. Then $\Cal U_{t}^{1}$ contains infinitely many 
smooth rational curves by $\vert \text{Aut}(\Cal U_{t}^{1}) 
\vert = \infty$ and by [PSS]. (See also [St] and [Kv].) Now the induced family 
over a generic linear section $0 \in \Delta \subset \Cal K^{1}$ satisfies the 
required property of (1). \qed \enddemo 

\demo{Proof of (2)} Let us take the subset 
$0 \in \Cal K^{2} \subset \Cal K = \Cal P$ defined by the equations 
$$(e + 2f, \omega) = (v_{11} + v_{13}, \omega) = (v_{21} + v_{23}, \omega) 
= 0.$$ 
Consider the induced family $u^{2} : \Cal U^{2} \rightarrow \Cal K^{2}$. Here 
$\Cal K^{2}$ is of dimension $17$. By construction, $\text{NS}(\Cal U_{t}^{2}) 
\simeq \langle e + 2f, v_{11} + v_{13}, v_{21} + v_{23} \rangle$ and the 
intersection matrix is: 
$$\pmatrix  4 & 0 & 0\\
            0 & -4 & 0\\
            0 & 0 & -4\\ 
\endpmatrix ,$$
if $t \in \Cal K^{2}$ is generic. 
Since the quadratic form 
$$Q_{2}(x, y) := 4x^{2} - 4y^{2} - 4z^{2}$$ 
does not represent $-2$, $\Cal U_{t}^{2}$ contains no smooth rational curve. 
In particular, the ample 
cone coincides with the positive cone. This already implies 
$\vert \text{Aut}(\Cal U_{t}^{2}) \vert = \infty$, because the automorphism group is then isomorphic to the orthogonal group of the Picard lattice up to finite groups [PSS]. Moreover, since 
$Q_{2}(x, y)$ admits infinitely many primitive integral zero's, $\Cal U_{t}^{2}$ admits infinitely many elliptic pencils. Now the induced 
family over a generic linear section $0 \in \Delta \subset \Cal K^{2}$ 
satisfies the required property of (2). \qed \enddemo 

\demo{Proof of (3)} Let us take the subset $\Cal K_{3} \subset \Cal K = 
\Cal P$ defined by the equations 
$$(e, \omega) = (f, \omega) = (v_{11} + v_{13} + v_{15} + v_{17}, \omega) 
= 0.$$ 
Consider the induced family $u^{3} : \Cal U^{3} \rightarrow \Cal K^{3}$. Here 
$\Cal K^{3}$ is of dimension $17$. By construction, one has 
$\text{NS}(\Cal U_{t}^{3}) \simeq U \oplus \langle -8 \rangle$ 
if $t \in \Cal K^{3}$ is generic. Then $X$ admits a Jacobian pencil 
$f : X \rightarrow \bold P^{1}$, i.e. an elliptic pencil with section, 
given by 
$U$. Here we set $X := \Cal U_{t}^{3}$. Since $f$ has no reducible fibers, the 
Mordell-Weil lattice $M(f)$ of $f$ is isomorphic to the lattice 
$\langle 8 \rangle$. (See the excellent survey of Shioda [Sh] for basic properties of the Mordell-Weil lattice.) Let $g$ be the generator of $M(f)$ corresponding to $8$. Then $g$ 
defines an automorphsim of $X$ of infinite order. In particular, 
$\vert \text{Aut}(X) \vert = \infty$. Let $C_{0}$ be the zero-section of $f$ and 
let $C_{1}$ be the section corresponding to $g$. Then $C_{1} = g(C_{0})$. 
Moreover, $g^{n}(C_{0})$ are mutually distinct smooth rational curves on $X$. 
On the 
other hand, by the definition of the height pairing $\langle * , ** \rangle$, 
we calculate 
$$8 = \langle C_{1}, C_{1} \rangle = 4 + 2(C_{0},C_{1}).$$ 
This gives $(C_{0},C_{1}) = 2$ and $(C_{n},C_{n+1}) = 2$. Thus, 
$C_{n} + C_{n+1}$ is a divisor of Kodaira type $I_{2}$ or $III$ and gives 
an elliptic pencil 
$$f_{n} := \Phi_{\vert C_{n} + C_{n+1} \vert} : X \rightarrow \bold P^{1}.$$ 
Since an elliptic pencil on a K3 surface admits at most 24 singular fibers, 
one has $\vert \{f_{n}\}_{n \geq 1} \vert = \infty$. 
Combining all these together, one finds that $\Cal U_{t}^{3}$ has 
$\rho(\Cal U_{t}^{3}) = 3$, infinitely many smooth rational curves, infinitely 
many elliptic pencils, and satisfies $\vert \text{Aut}(\Cal U_{t}^{3}) \vert 
= \infty$ if $t \in \Cal K^{3}$ is generic. Now the induced family over a 
generic linear section $0 \in \Delta \subset \Cal K^{3}$ satisfies the 
required property of (3). \qed \enddemo 
 
\demo{Proof of (4)} Let us take the subset $\Cal K^{4} \subset \Cal K = 
\Cal P$ defined by the equations 
$$(e + 6f, \omega) = (v_{11} + v_{13}, \omega) = (v_{21} + v_{23}, \omega) 
= 0.$$ 
Consider the induced family $u^{4} : \Cal U^{4} \rightarrow \Cal K^{4}$. Here 
$\Cal K^{4}$ is of dimension $17$. By construction, 
$\text{NS}(\Cal U_{t}^{4}) \simeq \langle e + 6f, v_{11} + v_{13}, 
v_{21} + v_{23} \rangle$ and the intersection matrix is: 
$$\pmatrix  12 & 0 & 0\\
            0 & -4 & 0\\
            0 & 0 & -4\\ 
\endpmatrix ,$$
if $t \in \Cal K^{4}$ is generic.  
Since the quadratic form 
$$Q_{4}(x, y) := 12x^{2} - 4y^{2} - 4z^{2}$$ 
represents neither $0$ nor $-2$, $\Cal U_{t}^{4}$ has neither smooth rational 
curve nor elliptic 
pencil. Then as before, the ample cone coincides with the positive cone and 
one 
has $\vert \text{Aut}(\Cal U_{t}^{4}) \vert = \infty$. Now the induced family over a generic linear section $0 \in \Delta 
\subset \Cal K^{4}$ satisfies the required property of (4). \qed \enddemo 

\demo{Proof of (5)} Let us take the subset $\Cal K^{5} \subset \Cal K = 
\Cal P$ defined by the equations 
$$(e, \omega) = (f, \omega) = (v, \omega) = 0.$$ 
Consider the induced family $u^{5} : \Cal U^{5} \rightarrow \Cal K^{5}$. Here 
$\Cal K^{5}$ is of dimension $17$ and one has $\text{NS}(\Cal U_{t}^{5}) = 
U \oplus A_{1}(-1)$ for generic $t \in \Cal K^{5}$ and $\Cal U_{t}^{5}$ enjoys 
all the properties required in (5) by [Ni3]. Now the induced 
family over a generic linear section $0 \in \Delta \subset \Cal K^{5}$ 
satisfies the required property of (5). Note that no element $h \in U \oplus 
A_{1}(-1)$ is ample on $S = \Cal U_{0}^{5}$, because such $h$ is perpendicular 
to $E_{8}(-1) \oplus E_{8}(-1)$. Therefore, the family is not projective. \qed 
\enddemo 
Now we are done. Q.E.D. for Proposition. \qed  
\enddocument